\documentclass[11pt,twoside,english]{article}
\usepackage{babel,amsfonts}
\usepackage{graphicx,epsfig,psfrag}

\newtheorem{Theorem}{Theorem}
\newtheorem{Proposition}{Proposition}[section]

\newtheorem{Definition}[Proposition]{Definition}

\newtheorem{Remark} [Proposition]{Remark}
\newtheorem{Oldtheorem}[Proposition]{Theorem}
\newtheorem{Paragraph}[Proposition]{}
\newcommand{\A}[7]{\bibitem[#7]{#1}{#2: \em #3. \em \small
    #4 (#5) \normalsize #6}}

\newcommand{\finalprueba}{ \begin{flushright}
  \rule{2mm}{2mm}
\end{flushright}
\vskip .5cm }

\begin{document}

\title{The  real analytic Feigenbaum-Coullet-Tresser attractor
 in the disk.}
\author{Eleonora Catsigeras\thanks{Corresponding author: eleonora@fing.edu.uy} , Marcelo Cerminara and Heber Enrich
   }

\date{}
\maketitle

\begin{center}

Instituto de Matem{\'a}tica (IMERL), Facultad de Ingener{\'\i}a

Universidad de la Rep{\'u}blica.

P.O.Box: C.C.30 Montevideo, Uruguay.

Tel-Fax +598-2-711-0621

E-mail addresses:

eleonora@fing.edu.uy

 cerminar@fing.edu.uy

  enrich@fing.edu.uy

\vspace{1cm}

\small
  May 2nd., 2008

\normalsize

\end{center}

\begin{abstract}
We consider a real analytic diffeomorphism $\psi _0$ on a $n$-dimensional disk ${\cal D}$, $n \geq 2$, exhibiting  a
Feigenbaum-Coullet-Tr{\'e}sser (F.C.T.) attractor, being far, in the $C^{\omega}({\cal D})$ topology,
 from the standard F.C.T.  map $\phi_0$  fixed by the double
renormalization.

We prove that  $\psi_0$ persists along
a  codimension-one
manifold ${\cal M} \subset C^{\omega}({\cal D})$, and that
 it is the
 bifurcating map along any one-parameter family  in
$C^{\omega}({\cal D})$ transversal to ${\cal M}$, from diffeomorphisms  attracted to
sinks, to those which exhibit chaos.

 The main tool in the proofs is a  theorem of Functional Analysis, which
 we state and prove in this paper, characterizing the existence of codimension
 one submanifolds in any abstract functional Banach space.

\end{abstract}



\section{Introduction}

 In dimension one, the Feigenbaum-Coullet-Tresser (F.C.T.) theory
\cite{CT, F, Fe} states that the F.C.T. attractor is a codimension
one phenomenon when seen locally, in a small neighborhood of the
real analytic standard map $\varphi _0$. This standard F.C.T. map $\varphi_0$, is the real
analytic unimodal map in the interval, of quadratic type at its
critical point, and such that $\varphi_0$ is fixed by the doubling
renormalization in the interval.

The existence of a local codimension-one manifold through
$\varphi_0$ is a consequence of the hyperbolicity of $\varphi _0$
as  fixed by the doubling renormalization in the space of real
analytic maps of the interval. (See the proofs in \cite{La, Su, Ly}).
This hyperbolic behavior was first proved by Lanford  in \cite{La}.

The codimension one character is also true for $n$-dimensional
real analytic transformations, as proved by Collet-Eckman and Koch
in \cite{CEK}, taking as the fixed point of the doubling
renormalization the standard endomorphic F.C.T. map $\phi_0: {\cal D} \mapsto \mbox {int
}({\cal D})$. This endomorphism  is defined from the standard F.C.T. map $\varphi_0$ in the interval, by endowing the
$2 \leq n$-dimensional disk ${\cal D}$ onto its interior, with infinite
codimension one contraction,  such that $\phi_0({\cal D})$ is
the graph of the map $\varphi_0$ in the interval. The precise
definition of the map $\phi_0$ will be reviewed in Definition
\ref{stardardFeigenMap} of this paper.

\vspace{.3cm}

Locally, nearby the standard F.C.T. map $\varphi _0$ in the interval, or nearby
the standard F.C.T. endomorphism $
\phi_0$ in the $n$-dimensional disk ${\cal D}$, the codimension one character
of infinitely doubling renormalizable maps was proved in the
space of $C^r$ transformations, provided that $r$ is large enough.
(See \cite{Da}, \cite{Su} and \cite{CE}).

Also for $C^r$ maps, far away  from the standard endormorphism $\phi _0$ in the
$n$-dimensional disk, if $n \geq 2$ and $r \geq 8$, the F.C.T.
attractor is a codimension one phenomenon, as proved in
\cite{CE99}). This result is not true in dimension one.  For
global results for maps in the interval exhibiting a F.C.T.
attractor see \cite{Su} and \cite{Ly}).

\vspace{.3cm}

We address here to the remaining open question about the codimension one character of the F.C.T.
attractor in the space $C^{\omega }({\cal D})$,  exhibited by infinitely doubling renormalizable  diffeomorhphisms or endomorphisms
$\psi _0$, that are far away from the standard F.C.T. endomorphism $\phi_0$ fixed by the renormalization. We prove that, in fact, it is  a codimension-
one phenomenon.

The key condition to obtain that result,   is the dimension two or greater of the manifold.
 The result is not true in dimension one.  In fact, the
known proof of the codimension one character of the F.C.T.
attractor in the $n-$disk ${\cal D}$ (see \cite{CE99}), requires the use of at
least two different spatial directions: a first direction onto which
asymptotically the map contracts after successive many doubling
renormalizations; and a second direction to perturb the map and to
construct a one-parameter family destroying the F.C.T. attractor.
This argument was used to prove that this attractor is a one-parameter bifurcating
phenomenon in $C^r({\cal D})$ in \cite{CE99}, and we will base on it our new result in $C^{\omega}({\cal D})$. The one parameter family such constructed shall be
transversal to the desired codimension-one manifold of maps
exhibiting the continuation of the F.C.T. attractor.

Nevertheless the known proof of the codimension-one character of
the F.C.T. attractor in the space  $C^r({\cal D})$,
 does not work for real analytic maps. This is due to
the construction, used in \cite{CE99}, of a one-parameter family of non-zero
diffeomorphisms or endomorphisms that have null derivative in infinitely many points of the compact disk ${\cal D}$.
(see Lemma 5.3 of the proof of Theorem 2 in \cite{CE99}).

Following the remark of Tresser (\cite{tresser00}),
and using the result of Theorem 2 in \cite {CE99} combined with
the density
 of the real analytic maps in the space  $C^r({\cal D})$ maps, we obtain the
 following Theorem, whose proof is the first purpose of this paper:

\begin{Theorem} \label{teoremaCodimension}
If $\psi _0: {\cal D} \mapsto \mbox {int}({\cal D})$ is a real
analytic map of the $n-$dimensional compact disk ${\cal D}$ to its
interior, where $n \geq 2$, and if $\psi _0$  exhibits a F.C.T.
attractor (see Definition \ref{FeigenAttractor}), then there exist
a local codimension-one $C^1$ manifold ${\cal M }$ in the Banach
space $C^{\omega} ({\cal D})$ of real analytic maps from ${\cal
D}$ to $\mathbb{R}^n$, such that $\psi _0 \in {\cal M}$ and for all $\psi \in {\cal M}$
the map $\psi :{\cal D} \mapsto \mbox {int}({\cal D})$ also
exhibits a F.C.T. attractor.
\end{Theorem}

We prove  Theorem \ref{teoremaCodimension}  in section \ref{pruebaTeoremas}.

\vspace{.4cm}

We base our arguments  on a  result stated in Theorem \ref{lemma} of this paper, which applies the   classical tools on Functional Analysis, dealing with non linear infinite dimensional submanifolds of an abstract Banach space $H$. More precisely, in Theorem \ref{lemma} we give necessary and sufficient conditions to obtain codimension-one submanifolds of a functional Banach space $H$, where some phenomenon appears, in terms of the persitence of the bifurcating quality of this phenomenon, along  one-parameter families in $H$.

The same arguments also work to prove  the following second result:

\begin{Theorem} \label{teoremaBifurcacion}
If $\psi_0$ verifies the hypothesis of Theorem
\ref{teoremaCodimension}  then $\psi _0$ belongs to a
one-parameter family of real analytic maps in ${\cal D}$ such that
at $\psi _0$ there exists a global bifurcation from maps that pass
through a cascade of period doubling bifurcations to maps that
pass through a sequence of homoclinic tangencies bifurcations.

Even more,
any one-parameter family that is transversal at $\psi_0$ to the manifold ${\cal M}
\subset C^{\omega}({\cal D})$ of the thesis of Theorem
 \ref{teoremaCodimension}, has the property above.
\end{Theorem}

We prove Theorem \ref{teoremaBifurcacion} in the last part of Section \ref{pruebaTeoremas}.

\vspace{.4cm}

The result of Theorem \ref{teoremaBifurcacion}, restricted to one-dimensional quadratic unimodal maps in the interval that are near  the standard F.C.T. map $\varphi_0$ (see Definition \ref{definicionStandardIntervalo}), was first obtained in \cite{EW}. Afterwards, the result was generalized to the disk in dimension $n \geq 2$ in \cite{C}, but also, only for maps that are in a small neighborhood of the standard F.C.T. map $\phi_0$ (see Definition \ref{stardardFeigenMap}).

As a consequence of Theorem \ref{teoremaBifurcacion},
 any  $\psi _0 $ showing a F.C.T. attractor, even far from the standard
 F.C.T. attractor, is the
 bifurcating map along any one-parameter family  in
$C^{\omega}({\cal D})$ transversal to ${\cal M}$.
At one side of $\psi_0$ the maps of the family exhibit
sinks, in a cascade of period doubling bifurcations, while
 at the other side they exhibit chaos (hyperbolic horseshoes, and also H{\'e}non-like attractors), due
to the sequence of homoclinic bifurcations that accumulate on $\phi_0$ (see \cite{PT}).

The conclusion of the F.C.T as a generic route to chaos,  is widely known and  applied in other sciences to physical autonomous dynamical
systems, but no mathematical proof
of it was known before, in open sets of $C^\omega ({\cal D})$, far away from
 the standard F.C.T. map.

 \vspace{.2cm}

 Let us  suggest that a similar result to that in the thesis of Theorem  \ref{teoremaCodimension}, can be obtained for other kind of infinitely renormalizable diffeomorphisms or endomorphisms in $n \geq 2$ dimensions. In fact,  instead of considering the classical standard F.C.T. map,  we can look at other fixed map by the doubling renormalization in the interval, with a non quadratic critical point or  in the critical circle (\cite{Yampolsky02}, \cite{Yampolsky06}). We would define other classes of Cantor set attractors that are not the F.C.T. attractor.

 Provided that  the  unidimensional map fixed by the renormalization also has a local hyperbolic behavior in the functional space (\cite{Yampolsky02}, \cite{deMdeF}), it defines  the corresponding endomorphisms in $n \geq 2$ dimensions, fixed and locally hyperbolic by the renormalization, with the same arguments used in \cite{CEK}.

 Finally, the technique tools we use in this paper can be applied to a diffeomorphism $\psi$ in $n$ dimensions, exhibiting the Cantor set attractor, but which is initially far away from that fixed and locally hyperbolic endomorphism.  We shall assume that its sequence of renormalized maps converges to that fixed endomorphism, in the functional space. We note that the renormalization is neither a linear operator nor a Fr{\'e}chet differentiable transformation in the functional space. Nevertheless, to construct a manifold ${\cal M}$ of the thesis of Theorem \ref{teoremaCodimension}, where the infinitely renormalizable Cantor set attractor persists, the main arguments in \cite{CE99} and in this paper can be applied. We are aware that the theory does not hold in dimension one.

\vspace{0.3cm}

\section{Definitions and previous results.}

A \em analytic $n$-disk ${\cal D}$, \em (or simply a \em disk
\em), is
 the image by a real analytic diffeomorphism of the unit
closed ball of $R^n$. ($n \geq 2$). Note that we call \em real
analytic diffeomorphism  \em to a real analytic transformation
that is invertible and whose inverse is also a real analytic
transformation.

Analogously a \em $C^r-n-$disk \em is the image  by a $C^r$
diffeomorphism of the unit closed ball of $R ^ n$. All analytic
disk is a $C^r$-disk.

\begin{Paragraph} \em {\bf The standard Feigenbaum-Coullet-Tresser map in the interval.} \label{definicionStandardIntervalo}

We call  \em  the standard F.C.T. map in the interval, \em to the unique
real analytic unimodal map $\varphi_0: [-1,1] \mapsto [-1,1]$ such
that $\varphi_0(0)=1, \; \varphi'(0)=0, \; \varphi_0''(0) < 0 $
and
$$\varphi_0(1)^{-1} \cdot \varphi_0 \circ \varphi_0 (\varphi_0(1) \cdot x) =
\varphi_0 (x) \; \; \; \forall x \in [-1,1] $$
\end{Paragraph}

The existence and unicity of $\varphi_0$ was the central
conjecture of the F.C.T. theory \cite {CT, F, Fe} and was proved
in \cite {La, Su}. We denote $\lambda $ to the number
$-\varphi_0(1) = 0.3995\ldots $.  The map $\varphi_0$ has a single
fixed point in $[-1,1]$, which  is larger than $\lambda $. The
analytic map $\varphi_0$ is symmetric: $\varphi_0 (x) = g_0 (x^2)
$ where $g_0$ is an analytic diffeomorphism from $[0,1]$ to
$[-\lambda,1]$. It can be analytically uniquely extended to an
open interval.

There exists one single periodic orbit of $\varphi_0$ of period
$2^N$ for each natural $N \geq 0 $, and this orbit is a hyperbolic
repeller. The orbit of countably many points eventually fall on
one of these repellers. All the other orbits of $\varphi_0 $ are
attracted to a Cantor set $K$ in the interval which we call
 \em the standard F.C.T. attractor in the interval. \em

\vspace{.2cm}

 Let $n \geq 2$. Let ${\cal D}$ be n-dimensional compact
 disk containing the segment $[-\lambda, 1] \times \{0 \}^{n-2} \times
[-\lambda, 1]$.

\begin{Definition}
\em  \label{stardardFeigenMap}{\bf The standard Feigenbaum-Coullet-Tresser
map in $n$-dimensions.}

 The map  $\phi_0 :{\cal D} \mapsto
\mbox{int}({\cal D})$ defined as
$$\phi _0 (x_1, x_2, \ldots, x_{n-1}, x_n ) = (x_n, 0, \ldots, 0,
\varphi_0(x_n))= (x_n, 0, \ldots, 0, g_0(x_n ^2) )$$ for all $x
\in {\cal D}$, will be called  \em the standard F.C.T. map in $n$
dimensions. \em It inherits the Cantor set attractor of the map
$\varphi_0$, that we call the \em standard F.C.T. attractor in $n$
dimensions.
\end{Definition}

\begin{Remark} \em \label{Rem9}
Observe that the standard F.C.T. map in $n$ dimensions has a one
dimensional character: it is an endomorphism of  ${\cal D}$
endowing  it to a one-dimensional image, contained in its
interior, and following the graph of $\varphi_0$.

The repellers of $\varphi_0$ are transformed into periodic
hyperbolic saddles of $\phi _0 $ with infinite contraction along
their stable manifolds. There
 exist
such a periodic orbit with period $2^N$ for each natural $N \geq 0
$. The unstable manifolds of the saddles have dimension one and
are contained in $\phi _0 ({\cal D})$. The stable manifold of each
saddle is the union of their pre-images by $\phi _0 $, formed by the intersection
 with ${\cal
D}$, of the
horizontal  $(n-1)$-dimensional hyperplanes for which $x_n$ is constant.

All the orbits of $\phi _0 $, except those in the stable manifolds
of the saddles, are attracted to its standard F.C.T. attractor.

\em \end{Remark}

We are interested in studying some Cantor set attractors for other
$n$-dimensional maps,  particularly for diffeomorphisms that might
be far away from the standard F.C.T. map.

\begin{Paragraph} \em {\bf Functional spaces.}

Given a analytic $n$-disk ${\cal D}$, the space $C^{\omega }({\cal
D})$ is the open set of the Banach space of  the real analytic
maps $\psi: {\cal D} \mapsto R^n$, such that $\psi ({\cal D})
\subset \mbox{int} ({\cal D})$. The   topology in
$C^{\omega}({\cal D})$ is that given by the supreme norm $\| \;
\psi \; \| = max \{\|\psi (x ) \|_{R^n}: x \in {\cal D}\}$.

Analogously the space $C^r({\cal D})$ is the open set of the
Banach space of all the $C^r $ maps $\psi: {\cal D} \mapsto R ^n$,
such that $\psi ({\cal D}) \subset \mbox{int} ({\cal D})$. The
topology in $C^{r}({\cal D})$ is that given by the supreme norm
$\| \; \psi \; \|_r = max \{\|\psi (x ) \|_{R^n}, \|D\psi (x )
\|, \|D^2\psi (x ) \|,\ldots, \|D^r \psi (x ) \| : x \in {\cal
D}\}$.

\em
\end{Paragraph}

In some parts of this paper we will need to work with the whole
Banach space of real analytic maps, or of $C^r$ maps, from ${\cal
D}$ to $R^n$ although their images are not contained in the
interior of ${\cal D}$. We will still denote them as
$C^{\omega}({\cal D})$ and $C^{r}({\cal D})$ , if there were no
risk of confusion.

\begin{Definition} \em {\bf Doubling renormalization.}
  \label{doubl}

  A map $\psi \in C^{\omega} ({\cal D})$ ($\psi \in C^{r} ({\cal D})$) is \em
doubling renormalizable  \em if there exists a analytic (resp.
$C^r$) $n-$disk ${\cal D}_1 \subset \mbox{int} {\cal D}$ such
that:
\begin {description}
\item
$\psi({\cal D}_1) \cap {\cal D}_1 = \emptyset $
\item
$\psi ^2 ({\cal D}_1) \subset \mbox{int}  ({\cal D}_1)$
\end{description}

If $\psi $ is doubling renormalizable and $\xi :{\cal D} \mapsto {\cal D}_1
 $ is a
real analytic (resp $C^r$) diffeomorphism (called \em change of
variables\em), the map $\cal{R} \psi $ defined as
      $ {\cal{R}} \psi = \xi ^{-1} \circ \psi \circ \psi \circ \xi $
is a \em renormalized map \em of $\psi $. \em
\end{Definition}

Note that doubling renormalizability is an open condition in
$C^{\omega} ({\cal D}) $ (resp. $C^r({\cal D})$). Also note that
${\cal R} \psi $ is not uniquely defined: small perturbations of
the change of variables $\xi $ give other renormalized map of
$\psi $. When referring to the properties of ${\cal R} \psi $ we
understand that there exists some renormalized map of $\psi $
having these properties.

\vspace{0.5cm}

By induction we define:

 {\bf $m$-times doubling renormalizable maps.}

A map $\psi \in C^{\omega } ({\cal D}) $ (resp. $\psi \in C^{r}
({\cal D}) $) is \em $m$-times (doubling) renormalizable \em if it
is $m-1$-times (doubling) renormalizable and its
$m-1$-renormalized ${\cal{R}}^{m-1} \psi $ is doubling
renormalizable. It is defined a \em $m$-renormalized map of $\psi
$ \em as ${\cal{R}}^m \psi = {\cal{R}}{ \cal {R}}^{m-1} \psi $

\vspace{.2cm}

 {\bf Infinitely doubling renormalizable maps.}

A map $\psi \in C^{\omega } ({\cal D}) $ (resp. $\psi \in C^{r}
({\cal D}) $) is \em infinitely (doubling)
renormalizable \em if it
is $m$-times (doubling) renormalizable for all natural $m $.

\vspace{.2cm}

The main example of infinite doubling renormalizable maps in the $n$-disk is the
fixed map standard F.C.T. endormorphism $\phi _0 $, defined in \ref{stardardFeigenMap}.

\begin{Remark} \em
An infinitely (doubling) renormalizable  diffeomorphism or endormorphism $\psi \in C^{\omega}({\cal D})$ (resp. $ C^{r}({\cal D})$) in the $n-$disk ${\cal D}$, may in general be far away from the stardard F.C.T. endormorphism $\phi_0$.
Nevertheless its sequence, or some subsequence, of some renormalized maps ${\cal R}^n \psi$ may converge to $\phi_0$ in the $C^{\omega}({\cal D})$ (resp. $C^r({\cal D})$) topology. If this happens, roughly speaking  $\psi$ exhibits a Cantor set attractor, that assympthotically inside, looks like the standard F.C.T. attractor. In more precise words, this is the statement of the following theorem:
\end{Remark}

\begin{Oldtheorem} \label{Prop1}
If a map  $\psi \in C^{\omega}({\cal D})$ (resp. $\psi \in
C^{r}({\cal D})$) is infinitely doubling renormalizable and if
there is a sequence ${\cal R}^j \psi$ of $j-$times
renormalizations of $\psi$ such that $$\lim _{j
\rightarrow \infty} {\cal R}^j \psi = \phi_0$$ in $C^{\omega }({\cal
D})$ (resp. in $ C^{r}({\cal D})$), where $\phi _0$ is the F.C.T.
map in $n-$ dimensions, then there exists a minimal Cantor set $K
\in {\mbox{int}(\cal D)}$ such that $\psi (K) = K$, there exits a
nighborhood $U \subset {\mbox{int}(\cal D)}$ of $K$ such that $K$
attracts almost all the orbits in $U$ and there exists a single
periodic saddle orbit in $U$ of period $2^N$ for all natural
number $N$ large enough.
\end{Oldtheorem}

{\em Proof:} See Theorem 2.12 of \cite{CE99} for the proof in $C^r({\cal D})$.
Exactly the same arguments of the proof in $C^r({\cal R})$ work in $C^{\omega}({\cal D})$, due that the topology of $C^{\omega}({\cal D})$ is  induced  from the $C^r({\cal D})$ topology.

\vspace{.5cm}

\begin{Definition} \em \label{FeigenAttractor} {\bf The F.C.T.
attractor.}

A map $\psi \in C^{\omega}({\cal D})$ (resp. $\psi \in C^{r}({\cal
D})$) from the $n-$ dimensional disk ${\cal D}$ to its interior
has a \em real analytic F.C.T. attractor \em if $\psi$ is
infinitely doubling renormalizable and there is a sequence ${\cal
R}^j \psi$ of $j-$times renormalizations of $\psi$, such that
$${\cal R}^j \psi \;\; _{j \rightarrow \infty} \rightarrow
 \phi_0$$ in $C^{\omega }({\cal D})$ (resp. in $ C^{r}({\cal D})$), where
$\phi _0$ is the standard F.C.T. map in $n-$ dimensions.

\end{Definition}

\begin{Remark} \em
 The Cantor set attractor $K$ of Theorem \ref{Prop1} and
 Definition \ref{FeigenAttractor}
 has bounded geometry in the sense that the diameter of the
connected compact atoms that asymptotically define $K$, decrease with an
asymptotic
rate  below 1. When looking microscopically the
decreasing rate, it tends to the number $\lambda = 0.3995\ldots $,
that is a spatial universal constant defined for the standard F.C.T map
$\phi_0$. In fact, $\lambda $ is the contraction rate of the
change of variables to renormalize $\phi_0$, and it is also the
asymptotic contraction rate of the change of variables to pass
from the ${\cal R}^j \psi$ to ${\cal R}^{j+1} \psi $, if $\psi$
verifies the Definition \ref{FeigenAttractor}. \em
\end{Remark}

We note  that Gambaudo and Tresser in \cite{GT} give an example
of a $n$ dimensional infinitely renormalizable map whose
renormalized maps do not converge to the F.C.T. map $\phi _0$. In
spite of that, this example has a Cantor set attractor that
verifies the thesis of the theorem \ref{Prop1}. Its geometry is
also bounded, but the bounds are different from $\lambda $. We do
not call that Cantor set a F.C.T. attractor.

\vspace{0.5 cm}

Let us recall some results from \cite{CE99} in which we will found
part of the proofs of Theorems \ref{teoremaCodimension} and
\ref{teoremaBifurcacion}:

\begin{Oldtheorem} \label{teoremac8}
For $r \geq 8$, if $\psi_0 \in C^r({\cal D})$ has a F.C.T.
attractor, then there exists a local codimension-one $C^1$ manifold
${\cal M}$ in $C^r({\cal D})$ such that $\psi _0 \in {\cal M}$ and
$\chi$ has a F.C.T. attractor for all $\chi \in {\cal M}$.
\end{Oldtheorem}

{\em Proof: } See Theorem 2 of \cite{CE99}.

\begin{Oldtheorem} \label{teoremabifurcacionc8}
For $r \geq 8$, if $\psi_0 \in C^r({\cal D})$ has a F.C.T.
attractor, then  any $C^1$ one-parameter family $\Psi = \{\psi
_t\}_{t \in [-1,1]} \subset C^r ({\cal D})$ transversal to the
manifold ${\cal M}$ of Theorem \ref{teoremac8} at $\psi _0$,
exhibits, at one side of $t=0$  a sequence of period doubling
bifurcations from sinks of period $2^N$ (for any sufficiently
large natural number $N$) to saddles of the same period and  sinks
of double period; and exhibits, at the other side of $t=0$,  a
sequence of homoclinic tangency bifurcations of saddles of period
$2^N$ (for any sufficiently large natural number $N$).
\end{Oldtheorem}

{\em Proof: } See Corollary 4 of \cite{CE99}.

\section {Characterization of codimension-one submanifolds in abstract Banach spaces.}

We recall  the following definitions, and the we state a theorem  dealing with the
differentiable submanifolds of codimension one, in abstract Banach spaces:

\begin{Definition} \label{topologyF}
\em {\bf (The Banach space of $C^1$ one-parameter families.)} Let
$H$ be a Banach space and let $\Psi = \{\psi _t\}_{t \in [-1,1]}$
be a $C^1$ one-parameter family in $H$, i.e. $\Psi $ is a $C^1$
application taking $t \in [-1,1] $ to $ \psi _ t \in H$. We denote
as $\partial \psi _t /\partial t \in H$ to the derivative respect
to $t $  of the application $t \mapsto \psi _t \in H$. The set $F=
C^1([-1,1], H)$ of all $C^1$ one-parameter families in $H$ is a
Banach space with the $C^1$ topology derived from the following
$C^1$ norm $$\|\Psi \|_F = max \{\| \psi _t
 \|_H, \|\partial \psi _t /\partial t  \|_H:  t \in [-1,1]\}$$
 We denote  $B_{\epsilon }(\Psi) $ to the open ball in $F$
 centered at $\Psi$ and with radius equal to $\epsilon >0 $.

 Given a $C^1$ one parameter family $\Psi = \{\psi_t\}_{t \in [-1,1]]} \in F $
 and given a
  fixed real number $t_0$ such that $|t_0| \leq 1$ we construct new
  families (many)
  in $ F $
   denoted as $(t_0)^*\Psi$ defined as:
   $$(t_0)^*\Psi = \{\widehat \psi _{t}\}_{t \in
    [-1,1]} \in F, \mbox{ where }   \widehat \psi _t = \psi _{t + t_0} \mbox{
    if }
    t \in [-1,1] , t + t_0 \in [-1,1]$$
     Note that, to define  $(t_0)^*\Psi = \{\widehat \psi _t\}_{t \in [-1,1]} = \{ \psi _{t + t_0}\}_{t \in [-1,1]} $, it is required to choose
       any $C^1$-extension
    of $\psi _{t+t_0}$, for the values of $t \in [-1,1]$ such that $t+t_0 \not \in
    [-1,1]$.

\end{Definition}

\begin{Definition} \label{definicionpersistencia}
\em {\bf (Persistent phenomena in $C^1$ one-parameter families.)}

Let $H$ be a Banach space, let ${\cal P}$ be any non empty subset
of $H$, and let $\Psi = \{\psi _t\}_{t \in [-1,1]} \in F$ be a
$C^1$ one-parameter family in $H$ such that $\psi _0 \in {\cal P
}$.

We say that \em the set $\cal P$ (or the phenomenon  ${\cal P}$)
is persistent in $C^1$ one-parameter families near $\Psi$ \em if
there exist $\epsilon
>0$ and a $C^1$ real function $a: B_{\epsilon }(\Psi) \subset F \mapsto
[-1 ,1 ]$ such that for all $\Gamma = \{\gamma _t\}_{t \in [-1,1]}
\in B_{\epsilon }(\Psi) \subset F$:

\begin{itemize}
\item [ a) ] $\gamma _{a(\Gamma )} \in {\cal P}$

 \item [ b) ] If
$\gamma _0 = \psi _0$ then $a (\Gamma) = 0$. (In particular $a
(\Psi) = 0$.)

\item [ c) ] If $|t_0| $ is small enough  then $a((t_0)^* \Gamma)
= a(\Gamma) - t_0$. (In particular $a((t_0)^* \Psi) = - t_0$.)
\end{itemize}
To explicit the value of $\epsilon$ in this definition we will refer to the
set ${\cal P}$ as being \em $\epsilon-$persistet in $C^1$ one-parameter families near
$\Psi$. \em

\end{Definition}

\begin{Oldtheorem}
 \label{lemma} Let $H$ be a Banach space, let ${\cal P}$ be any
non empty subset of $H$ and let  $\psi _0 \in {\cal P}$. The
following assertions are equivalent:

\begin{itemize} \em
\item [  i) ]  \em There exists a $C^1$ local manifold $\cal M$ in $H$ with
codimension one such that $\psi _ 0 \in {\cal M}$ and $\chi \in
{\cal P}$ for all $\chi \in {\cal M}$. \em
\item [ ii) ]  \em There exists a $C^1$ one-parameter family $\Psi = \{\psi _t\}
_{t \in [-1,1]}$ passing through $\psi _0$ for $t = 0$ and such
that
 the set $\cal P$ is  persistent in
$C^1$ one-parameter families near $\Psi$  \em (according with
definition \ref{definicionpersistencia}).
\item [iii) ] \em There exists $v_0 \in H$ such that the set ${\cal P}$
is persistent in $C^1$ one-parameter families near $\Psi =
\{\psi_0 + t v_0\}_{ t \in [-1,1]}$ \em (according with definition
\ref{definicionpersistencia}).

\end{itemize}

\end{Oldtheorem}

{\em Proof: }

We first prove  that i) implies ii):

 We apply the $C^1$
persistence of the transversal intersection between $C^1$ manifolds in H (see
\cite{Lang}): Being ${\cal M}$ a codimension one, $C^1$ manifold o $H$ passing through
$\psi_0$, it is locally characterized by a real equation $b(\chi) = 0$ in a
neighborhood of radius $\delta >0$ of $\psi _0 \in H$:
$${ \cal M} = \{ \chi \in H : \| \chi - \psi _0 \|_H < \delta, \;  b(\chi) =
0 \} $$  where $b: \{ \chi \in H : \| \chi - \psi _0 \|_H <
\delta \} \mapsto \mathbb{R}$ is some real function of $C^1$ class and with
surjective derivative.

Let $v_0 \in H$ be transversal to ${\cal M}$ at $\psi _0$. That is
$$Db_{\psi _0} \cdot v_0 \neq 0 $$

If $\|v_0 \| _H>0$ is small enough, then the family $\Psi = \{\psi_t\}_{t \in [-1,1] }= \{ \psi
_0 + t v_0\}_{t \in [-1,1]} \in F$ verifies $\|\psi _t -\psi
_0\|_H< \delta $ for all $t \in [-1,1]$.

Taking a smaller positive value for $\delta$, let us define the
transformation $G: B_{\delta } (\Psi) \times [-1 , 1 ] $ such that,
 if $\Gamma = \{\gamma_t\}_{t \in [-1,1]} \in B_{\delta }(\Psi
) \subset F$ and if $t \in [-1 ,1 ] $, then:
$$G(\Gamma , t ) = b(\gamma_t)$$

The transformation $G$ is $C^1$ because it is the composition of
the $C^1$ real function $b$ with the parameter evaluation
$\gamma_t$ of the $C^1$ parameter family $\Gamma$.

As $\psi _0 \in {\cal M}$ we have  $$G(\Psi, 0) = 0 , \; \;
\left . \frac {\partial G}{\partial t} \right |_{\Gamma = \Psi, t=
0} = Db |_{\psi _0} \cdot v_0 \neq 0
$$
Then, by the implicit function theorem there exists $\epsilon >0 $
and $a: B_{\epsilon } (\Psi ) \subset F \mapsto [-1 , 1 ]$ of $C^1
$ class, such that $$G(\Gamma, a(\Gamma )) = 0 $$ Thus
$b(\gamma_{a(\Gamma )}) =0$. Therefore $\gamma _{a (\Gamma )} \in
{\cal M} \subset {\cal {P}}$. We conclude that the $C^1$ real
function $a$ verifies condition a) of Definition
\ref{definicionpersistencia}.

The Implicit Function Theorem also asserts that if $G(\Gamma, a) =
0$ for some $\Gamma \in B_{\epsilon } (\Psi)$ and some $a$ with
$|a|$ small enough, then $a = a (\Gamma )$. This last assertion
proves that the real function $a$ verifies also
conditions b) and c) of Definition \ref{definicionpersistencia},
as wanted.

\vspace{1cm}

Let us now prove that ii) implies iii):

Let $\Psi = \{\psi _t\}_{ t \in [-1,1]} \in F$ be the $C^1$
one-parameter family given in ii) and $\epsilon >0$ the real
number given in the definition \ref{definicionpersistencia} of
persistence of ${\cal P}$ in one-parameter families near $\Psi$.
Let us call $v_0 = \left . ( \partial \psi / \partial t ) \right
|_{t=0}$. The $C^1$ condition of $\Psi$ implies that there exists
 $0 < \delta < 1$ such that for $ |t| \leq \delta $:
 $$\| \psi _0 + t v_0 - \psi _t \|_H < \epsilon /2 , \; \;
 \; \|\frac {\partial \psi _t}{\partial t} - v_0 \|_H < \epsilon /2 $$
Take any $C^1$ extension $\rho = \{\rho _t\}_ {t \in [-1,1]} \in
F$ of $\{\psi _0 + t v_0\} _{t \in [-\delta, \delta]}$ such that
$\rho \in B_{\epsilon/2}(\Psi) \subset F$.

As the phenomenon ${\cal P}$ is $\epsilon-$persistent in $C^1$ one
parameter families near $\Psi$, and $\rho$ is $\epsilon /2-$near
$\Psi$, we obtain that ${\cal P}$ is $\epsilon /2$ persistent in
$C^1$ one parameter families near $\rho$.

Now we shall  construct a family $\Lambda \in F $ to be linear on
the parameter $t \in [-1,1]$  as in the thesis iii), from $\rho
\in F$. (The one parameter family $\rho$ is linear only in the
small $\delta$- neighborhood of the parameter value $t= 0$.)

Consider $\Lambda = \{\lambda _t\}_{t \in [-1,1]}\subset F$
defined as $\lambda _t = \psi _0 + \delta t v_0$ for all $t \in
[-1,1]$.

As $\rho_0 = \psi _0$ then $a (\rho) = 0 $. By continuity of the
function $a$ there exists $\epsilon '$ such that $0 < \epsilon ' <
\epsilon /2$ and $|a (\tilde \Gamma)| < \delta /2 $ if $\| \tilde
\Gamma - \rho \| _ F < \epsilon ' $.

 Take $\epsilon '' = \epsilon ' \delta $. It is enough to
prove that the phenomenon ${\cal P}$ is $\epsilon '' $ persistent
in $C^1$ one-parameter families near $\Lambda$.

 In fact, for all $\Gamma = \{\gamma _t\}_{t \in [-1,1]} \in B_{\epsilon ''} (\Lambda)
 \subset F$ we have
 $$\|\gamma _t - (\psi _0 + \delta t v_0)\|_H < \epsilon '' = \epsilon ' \delta  <
 \epsilon '
 \; \; \; \; (1)$$
 $$\left \|\frac{\partial \gamma _t}{\partial t} - \delta v_0 \right \| _H
 < \epsilon '' = \epsilon' \delta , \; \; \; \; \left \| \frac{1}{\delta}\frac{
 \partial \gamma _t}{\partial t} - v_0 \right \|_H < \epsilon ' \; \; \; (2)$$
Observe that $\rho _t = \lambda _{t /\delta } = \psi _0 + t v_0$
for $|t|\leq \delta$. Analogously, for $\Gamma = \{\gamma _t\}_{t
\in [-1,1]} \in B_{\epsilon ''} (\Lambda) \subset F$ define
    $\tilde \gamma _t =  \gamma _{t /\delta }$ for $|t|\leq \delta $, and
 consider any $C^1$ extension $\tilde \Gamma = \{\tilde \gamma _t\}_{t \in [-1,1]} \in F$.

  Using inequalities (1) and  (2) we check that $\tilde
  \gamma_t
  $ is $C^1\epsilon'-$near $\rho$ for parameter values $|t| \leq
  \delta$, and so the extension $\tilde \Gamma \in F $ to the whole parameter
  domain
  $t \in [-1,1]$ could be chosen such that $\tilde \Gamma \in B_{\epsilon '}(\rho) \subset
  F$.

 The $\epsilon ' -$persistence of the phenomenon ${\cal
 P}$ in one-parameter families near $\rho $, for the chosen value of
 $\epsilon '< \epsilon /2$, allows the existence
 of a $C^1$ function $\tilde a: B_{\epsilon '} \mapsto [-\delta/2, \delta/2]$
verifying the conditions of the definition
\ref{definicionpersistencia}.

As the values of ${\tilde a }$ are contained in the
$\delta/2$-neighborhood of 0 in the parameter domain, we have that
$\tilde \gamma _{\tilde a (\tilde \Gamma)} = \gamma _{\tilde a
(\tilde \Gamma) /\delta } \in {\cal P}$. The map $\tilde a$
depends on the given $\Gamma \in B_{\epsilon '' }(\Lambda)$  and
is independent on the choice of the extension $\tilde \Gamma$ for
parameter values outside $[-\delta, \delta]$.

Define $a (\Gamma) = \tilde a (\tilde \Gamma) / \delta$. It is
straightforward to check that the real function $a$ verifies the
conditions of the definition \ref{definicionpersistencia}

\vspace{1cm}

Finally let us prove that iii) implies i):

Let $\Psi = \{\psi _ 0 + t v_0\}_{t \in [-1,1]} \in F$ be the
one-parameter family given in the hypothesis (iii). Let $\epsilon
>0$ be the radius of the ball centered at $\Psi$ in $F$, where the
$C^1$-real function $a$  is defined, according to  Definition
\ref{definicionpersistencia} of persistence of the phenomenon
${\cal P}$.

Let us choose $\delta >0$ such that if $\|\chi - \psi _0 \|_H <
\delta$ then $$\Gamma (\chi) = \{\chi + t v_0\}_{t \in [-1,1]} \in
B_{\epsilon }(\Psi) \subset F \; \; \; \; (3)$$

Let us define $b: B_{\delta }(\psi _0 ) \subset H \mapsto [-1,1]$
as $$ b(\chi) = a (\Gamma (\chi))$$ By construction we have
$b(\psi _0) = 0$ and
$$b(\chi ) = 0 \; \; \Rightarrow \; \; \chi \in {\cal P}$$

Our aim is to prove that the set ${\cal M } \subset H$, defined as
$${\cal M} = \{ \chi \in B_\delta (\psi _0) \subset H: b(\chi)= 0\},$$
is an embedded $C^1$ local manifold. It is immediate that $\psi _0
\in {\cal M }$. It is enough to prove that the real function $b$
is of $C^1$- class and that its derivative at $\chi = \psi_0$ (
$Db (\psi _0): H \mapsto R$)
 is surjective.

 First, the function $b $ is the composition $b (\chi) = a \circ \Gamma
 (\chi)$. As $a: B_{\epsilon } (\Psi) \subset F \mapsto R$ is of
 $C^1$ class by assumption, it is left to prove that  $\Gamma : B_{\delta } (\psi _0 )
 \subset H \mapsto B_{\epsilon }(\Psi) \subset F$ defined  in
 (3) is differentiable with continuous derivative.

 Let us take $\Delta \chi \in H$ such that $\chi + \Delta \chi \in B_{\delta }(\psi
 _0)$.
 From (3) the increment of the function $\Gamma $ in $\chi$ is

 $$\Gamma (\chi + \Delta \chi) - \Gamma (\chi) = \{\Delta \chi\}_{t \in [-1,1]} \in F$$

In other words,  the increment of $\Gamma $ is $\Delta \Gamma =
i(\Delta \chi)$, where $i$ is the inclusion defined as follows:
$$i: H \hookrightarrow F \; \mbox{ such that } \; \forall \gamma \in H:
i(\gamma) = \{\gamma_t \}_{t \in [-1,1]} \in F \; \; \mbox{ where
} \; \gamma _t = \gamma \; \forall t \in [-1,1] $$

As the inclusion $i: H \hookrightarrow F$ is linear and
continuous, then the map $\Gamma : B_{\delta } (\psi _0 )
 \subset H \mapsto B_{\epsilon }(\Psi) \subset F$ defined  in
 (3) is of $C^1$ class as wanted.

 Finally, it is left to prove that $Db (\Psi _0): H \mapsto \mathbb{R}$ is
 surjective. It is enough to prove that it is not null the directional derivative
 of the real function $b: B_{\delta }(\psi _0) \subset H \mapsto R$
  along the direction $v_0 \in H$ at $\psi _0$. (Note that although it was not
  asked $v_0 \neq 0 $ in the hypothesis
  iii), the assertion above also proves that iii) implies $v_0 \neq 0
  $.)

  In fact, for any
  real number $\lambda$ sufficiently small in absolute value so that $\psi _0 + \lambda v_0 \in
  B_{\delta } (\psi _0) \subset H$, we have:
  $$b(\psi_0 + \lambda v_0) = a(\{\psi _0 + \lambda v_0 + t v_0\}_{t \in [-1,1]})
  = a(\lambda ^* \Psi) = a(\Psi) - \lambda = - \lambda$$
Therefore
$$\frac{d}{d\lambda} b (\psi _0 + \lambda v_0) = \frac{d}{d\lambda} (-\lambda) = -1 \neq 0
 $$ \finalprueba

\section {Conclusion of the main results: Theorems \ref{teoremaCodimension} and \ref{teoremaBifurcacion}.}
\label{pruebaTeoremas}

 {\bf \em Proof of Theorem
\ref{teoremaCodimension}:} Let ${\cal P}$ be the set of
transformations exhibiting a F.C.T. attractor in $C^r({\cal
D})\subset C^{\omega}({\cal D})$. Let $\psi_0 \in {\cal P} \cap
C^{\omega}({\cal D})$.

Let $F^{r}$ be the space of $C^1$ one-parameter families of maps
in $C^r({\cal D})$ and let $F^{\omega}$ be the space of $C^1$
one-parameter families of maps in $C^{\omega }({\cal D})$. We have
that $F^{\omega} \subset F^{r}$ and the topology in $F^{\omega }$
defined in \ref{topologyF} is the induced topology from $F^{r}$.

Due to Theorem \ref{teoremac8} and Lemma \ref{lemma} there exists
$v_0 \in C^r({\cal D})$ such that the phenomenon ${\cal P}$ is
persistent in one-parameter families in $F^r$ near $\Psi =
\{\psi_0 + t v_0\}_ {t \in [-1,1]}$. Let $\epsilon >0 $ be the
number of the definition \ref{definicionpersistencia}.

As $C^{\omega}({\cal D})$ is dense in $C^r({\cal D})$ there exists
$w_0 \in C^{\omega }({\cal D })$ such that $\|w_0 -
v_0\|_{C^r({\cal D})} < \epsilon /2$. Define the one-parameter
family $\widetilde \Psi\in F^{\omega }$  as follows: $\widetilde
\Psi = \{\psi_0 + t w_0\}_{t \in [-1,1]} $.

We obtain now  that the phenomenon ${\cal P}$ is persistent in one
parameter families of $F^r$ near $\widetilde \Psi$, taking
$\epsilon/2$ instead of $\epsilon $ in the definition
\ref{definicionpersistencia}. Restricting to the subspace
$F^{\omega}$ we conclude that the phenomenon ${\cal P} \cap
C^{\omega}({\cal D})$ is persistent in one parameter families of
$F^{\omega}$ near $\widetilde \Psi$. Finally, applying Lemma
\ref{lemma} again, we obtain that there exists in
$C^{\omega}({\cal D})$ the local codimension one manifold ${\cal
M}$ of maps in ${\cal P}\cap C^{\omega}({\cal D})$ as wanted.
\finalprueba

{\bf \em Proof of Theorem \ref{teoremaBifurcacion}:} We
denote $F^r$ and $F^{\omega}$ as in the last proof. Applying
Theorem \ref{teoremabifurcacionc8} all one-parameter family $\Psi
\in F^r$ that is transversal to the local codimension one manifold
${\cal M}$ of maps exhibiting a F.C.T. attractor in $C^r({\cal
D})$ verifies the thesis of Theorem \ref{teoremaBifurcacion}. As
in the last proof, the density of $C^{\omega}({\cal D})$ in
$C^{r}({\cal D})$ implies that there exists a one-parameter family
$\widetilde \Psi \in F^{\omega}$ transversal to ${\cal M}$ in
$C^{r}({\cal D})$, and thus, verifying the thesis. \finalprueba

\end{document}